\documentclass[11pt]{amsart}
\usepackage{graphicx}
\usepackage{amssymb, amsmath}
\begin{document}
\title{A Note on Derivations of Lie Algebras}
\author{\sc M. Shahryari}
\thanks{{\scriptsize
\hskip -0.4 true cm MSC(2010): 17B40
}}

\address{ Department of Pure Mathematics,  Faculty of Mathematical
Sciences, University of Tabriz, Tabriz, Iran }
\email{mshahryari@tabrizu.ac.ir}
\date{\today}

\begin{abstract}
In this note, we will prove that a finite dimensional Lie algebra
$L$ of characteristic zero, admitting an abelian algebra of
derivations $D\leq Der(L)$ with the property
$$
L^n\subseteq \sum_{d\in D}d(L)
$$
for some $n>1$, is necessarily solvable. As a result, if $L$ has
a derivation $d:L\to L$, such that $L^n\subseteq d(L)$, for some
$n>1$, then $L$ is solvable.
\end{abstract}

\maketitle

In \cite{Lad}, F. Ladisch proved that a finite group $G$, admitting
an element $a$ with the property $G^{\prime}=[a, G]$, is solvable.
Using this result, one can prove that a finite group is solvable, if
it has a fixed point free automorphism. In this note, we prove a
similar result for Lie algebras in a more general framework; we show
that a finite dimensional Lie algebra $L$ of characteristic zero, is
solvable if it has an abelian subalgebra $A$ with the property
$L^n\subseteq [A, L]$, for some $n>1$. Next, we use this result
to prove that a finite dimensional Lie algebra $L$ of characteristic
zero, admitting an abelian algebra of derivations $D\leq Der(L)$
with the property
$$
L^n\subseteq \sum_{d\in D}d(L)
$$
for some $n>1$, is necessarily solvable. As a special case, we
conclude that if the Lie algebra $L$ admits a derivation $d:L\to L$,
such that $L^n\subseteq d(L)$, for some $n>1$, then $L$ is
solvable. Note that a similar result was obtained by N. Jacobson in
\cite{Jac}: a finite dimensional Lie algebra of characteristic zero, admitting an invertible derivation, is nilpotent.\\
Our main theorem (Theorem 1 bellow) is also true for connected
compact Lie groups and so, it may be also true for finite groups.
Therefore, we ask the
following question;\\

{\em Let $G$ be a finite group admitting an abelian subgroup $A$
with the property $G^n\subseteq [A, G]$, for some $n>1$. Is it
true that $G$ is solvable?}\\

During this note, $L$ is a finite dimensional Lie algebra over a
field $K$ of characteristic zero. By $L^n$ and $L^{(n)}$, we will
denote the $n$-th terms of the lower central series and derived
series of $L$, respectively. Also, $Der(L)$ will
denote the algebra of derivations of $L$.\\

{\bf Theorem 1.} Suppose there exists an abelian subalgebra  $A\leq L$ and an integer $n>1$,
such that $L^n\subseteq [A, L]$. Then $L$ is solvable.\\

{\bf Proof.} Let $S=L^{n-1}$. First, we show that $S$ is solvable.
To do this, we use Cartan criterion. Let $x\in S$ and $y\in
S^{\prime}$. Since $S^{\prime}=[L^{n-1},L^{n-1}]\subseteq
L^n\subseteq [A, L]$, so
$$
y=\sum_{i}[a_i, u_i],
$$
for some $a_1, \ldots, a_k\in A$ and $u_1, \ldots, u_k\in L$. Now,
we have
$$
Tr (ad_Sx\ ad_Sy)=\sum_i Tr (ad_Sx\ ad_S[a_i, u_i]).
$$
Since $S$ is an ideal, we have $ad_S[a_i, u_i]=ad_Sa_i\
ad_Su_i-ad_Su_i\ ad_Sa_i$ and hence
\begin{eqnarray*}
Tr (ad_Sx\ ad_Sy)&=&\sum_i Tr (ad_Sx\ ad_Sa_i\ ad_Su_i-ad_Sx\ ad_Su_i\ ad_Sa_i)\\
                 &=&\sum_i Tr (ad_Sa_i\ ad_S[u_i, x]).
\end{eqnarray*}
Now, $[u_i,x]\in [L, S]=L^n\subseteq [A, L]$, and so
$$
[u_i, x]=\sum_j [b_{ij},v_j],
$$
for some $b_{i1}, \ldots, b_{il}\in A$ and $v_1, \ldots, v_l\in L$.
Therefore
\begin{eqnarray*}
Tr (ad_Sa_i\ ad_S[u_i, x])&=&\sum_j Tr\ (ad_Sa_i\ ad_S[b_{ij}, v_j])\\
                          &=&\sum_j Tr\ (ad_S[a_i, b_{ij}]ad_Sv_j)\\
                          &=& 0.
\end{eqnarray*}
Therefore,
$$
Tr (ad_Sx\ ad_Sy)=0,
$$
and hence $S$ is solvable. We have $L^{(n-2)}\subseteq L^{n-1}$, so $L$ is solvable.\\

As a result, we have;\\

{\bf Corollary 1.} Suppose $L$ is semisimple and $A$ is an abelian
subalgebra. Then $[A, L]\subsetneqq L$.\\

Using Lie functor, we can restate Theorem 1, for
connected compact Lie groups;\\

{\bf Corollary 2.} Suppose a connected compact Lie group $G$ has an
abelian Lie subgroup $A$, such that $G^n\subseteq [A,G]$, for some
$n>1$. Then $G$ is solvable.\\

{\bf Theorem 2.} Suppose there is an abelian subalgebra $D\leq
Der(L)$ and an integer $n>1$ such that
$$
L^n\subseteq \sum_{d\in D}d(L).
$$
Then $L$ is solvable.\\

{\bf Proof.} Suppose $\hat{L}=D\ltimes L$. Note that, elements of
$\hat{L}$ are of the form $(d, x)$, with $d\in D$ and $x\in L$.
Also, we have
$$
[(d, x), (d^{\prime}, y)]=(0, [x, y]+d(y)-d^{\prime}(x)).
$$
It is easy to see that for $d_1, \ldots, d_n\in D$ and $x_1, \ldots,
x_n\in L$, we have
$$
[(d_1, x_1), \ldots, (d_n, x_n)]=(0, [x_1, \ldots, x_n]+y),
$$
for some $y\in \sum_{d\in D}d(L)$. Now,
$$
[x_1, \ldots, x_n]\in \sum_{d\in D}d(L),
$$
so there exists $\delta_1, \ldots, \delta_k\in D$ and $u_i, \ldots,
u_k\in L$, such that
$$
[x_1, \ldots, x_n]+y=\sum_i \delta_i(u_i).
$$
We have
\begin{eqnarray*}
[(d_1, x_1), \ldots, (d_n, x_n)]&=&(0, [x_1, \ldots, x_n]+y)\\
                &=&\sum_i (0, \delta_i(u_i))\\
                &=&\sum_i[(\delta_i, 0), (0, u_i)]\\
                &\in & [D, \hat{L}].
\end{eqnarray*}
Therefore, $\hat{L}^n\subseteq [ D, \hat{L}]$, and hence
$\hat{L}$ is solvable. So $L$ is also solvable.\\

As a special case, if the Lie algebra $L$ admits a derivation
$d:L\to L$, such that $L^n\subseteq d(L)$, for some $n>1$, then
$L$ is solvable.\\

{\bf Acknowledgment.} The author would like to thank P. Shumyatiski
and K. Ersoy for their comments and suggestions.


\end{document}